# ANALYSIS OF TOP-SWAP SHUFFLING FOR GENOME REARRANGEMENTS

By Nayantara Bhatnagar,[1] Pietro Caputo,[2] Prasad Tetali[3]
and Eric Vigoda[1]

*Georgia Institute of Technology, Università di Roma Tre, Georgia Institute of Technology and Georgia Institute of Technology*

We study Markov chains which model genome rearrangements. These models are useful for studying the equilibrium distribution of chromosomal lengths, and are used in methods for estimating genomic distances. The primary Markov chain studied in this paper is the top-swap Markov chain. The top-swap chain is a card-shuffling process with $n$ cards divided over $k$ decks, where the cards are ordered within each deck. A transition consists of choosing a random pair of cards, and if the cards lie in different decks, we cut each deck at the chosen card and exchange the tops of the two decks. We prove precise bounds on the relaxation time (inverse spectral gap) of the top-swap chain. In particular, we prove the relaxation time is $\Theta(n+k)$. This resolves an open question of Durrett.

**1. Introduction.** Genome rearrangements play an important role in a variety of biological studies, for example, genomic distance [8, 9], phylogenetic studies [11] and cancer biology [7]. Rearrangements refer to chromosomal fissions and fusions, reciprocal translocations between chromosomes and inversions within a chromosome. Fissions, fusions and reciprocal translocations are examples of translocations. In this paper we study Markov chains which model genomic rearrangements by translocations and inversions.

Stochastic models of chromosomal rearrangements by translocations were introduced by Sankoff and Feretti [10]. Such models are useful for studying the equilibrium distribution of chromosomal lengths. De, Durrett, Ferguson

Received September 2006; revised January 2007.
[1]Supported by NSF Grant CCF-0455666.
[2]Supported by MIUR PRIN–COFIN.
[3]Supported by NSF Grant DMS-04-01239.

and Sindi [3] studied refinements of the models, including the introduction of a fitness function on chromosomal lengths, and showed that these models have a reasonable fit with data for many species. These models were used by Durrett, Nielsen and York [6] in a Bayesian approach for estimating genomic distance.

Various of the above models being simply Markov chains, Durrett [5] considered estimating the rate of convergence to stationarity of the corresponding Markov chains. One of the open problems raised by Durrett pertains to the analysis of, what we call henceforth, the *top-swap* chain. We study the so-called *relaxation time* of the chain, which is the inverse of the spectral gap of the transition matrix of the chain. The relaxation time is the key quantity in the rate of convergence, and hence, it is of utmost importance in the efficiency of any simulations of the model. Beyond its computational significance, the convergence rate also has biological significance since it addresses the rate at which genomes reach equilibrium.

The top-swap chain has two parameters: the number of chromosomes, denoted as $k$, and the number of genes, denoted as $n$. The chain can be viewed as a card shuffling problem with $k$ decks. More precisely, the state space of the chain is a partition of $n$ distinct cards into $k$ decks. The cards have some ordering within each deck, and the decks are labeled. At each transition we choose two random positions, where if the size of the $i$th deck is $n_i$, there are $n_i + 1$ positions to choose from in the $i$th deck. If the positions are in the same deck we do nothing. Otherwise, we cut both decks at the chosen positions. We then exchange the tops of the two decks. The figure below illustrates a sample transition (for $k = 4, n = 12$) where the chosen pair of positions is marked:

$$
\begin{array}{cccccccc}
 &  &  &  &  &  & 4 &  \\
 & 9 &  &  &  &  & 3 &  \\
4 & 7 & \underline{2} &  &  & 7 & 5 &  \\
3 & 1 & 11 & 8 & \implies \ 9 & 1 & 11 & 8 \\
\underline{5} & 10 & 6 & 12 &  2 & 10 & 6 & 12.
\end{array}
$$

The Markov chain allows empty decks, hence, fusions and fissions are modeled by this stochastic process. Translocations with nonempty decks are known as reciprocal translocations. As we note below, inversions can also be included in the above Markov chain.

Let $\tau(n, k)$ denote the inverse of the spectral gap of the top-swap chain. This is known as the relaxation time. We prove tight bounds, up to constant factors, on the spectral gap of the top-swap Markov chain.

THEOREM 1.1. *There exist constants $C_1 > C_2 > 0$ such that*

$$C_1(n+k) \geq \tau(n,k) \geq C_2(n+k).$$

REMARK 1.1. Durrett [5] proposes a Markov chain which is identical to the top-swap chain, except when the chosen pair of positions is in the same deck, this interval of cards is inverted. The above upper bound on the relaxation time immediately applies to this chain, since comparison of Dirichlet forms shows that extra transitions can only decrease the relaxation time.

As a byproduct of our proof, we also obtain tight bounds on the relaxation time of a $k$-deck random transposition chain; this chain is a natural extension of the classical (1-deck) random transposition chain to $k$-decks, for $k > 1$. For a precise statement of the result, see Theorem 5.2.

Before giving a high-level description of our proof, we need to introduce some notation. Let $\Omega$ denote the state space of the top-swap chain, and $P$ denote its transition matrix. We let $\nu$ denote the uniform distribution over $\Omega$. Since $P$ is symmetric, $\nu$ is reversible. Therefore, $\nu$ is the stationary distribution of the chain. For any $f : \Omega \to \mathbb{R}$, the Dirichlet form of the process is defined as

$$\mathcal{E}(f) = \tfrac{1}{2} \sum_{\sigma,\eta \in \Omega} \nu(\sigma) P(\sigma,\eta)(f(\sigma) - f(\eta))^2$$

and the variance is

$$\operatorname{Var}(f) = \sum_{\sigma \in \Omega} \nu(\sigma)(f(\sigma) - \mathbf{E}_\nu(f))^2$$
$$= \tfrac{1}{2} \sum_{\sigma,\eta \in \Omega} \nu(\sigma)\nu(\eta)(f(\sigma) - f(\eta))^2.$$

We then have that the relaxation time is

$$\tau(n,k) = \sup_f \frac{\operatorname{Var}(f)}{\mathcal{E}(f)},$$

where the supremum is taken over all nonconstant functions $f$.

The lower bound on the relaxation time is easy to show by taking $f$ as the indicator function for whether the first deck is empty. We give a high-level sketch of the proof outline in Section 2 before presenting a detailed proof in Sections 3, 4 and 5.

**2. High-level proof description.** The analysis of the spectral gap of the top-swap chain has two major parts. The first part shows that it suffices to analyze the spectral gap for the 2-deck version of top-swap. (The 2-deck version of top-swap is simply the top-swap chain with $k = 2$.) We then analyze the 2-deck top-swap chain in the second part of our proof. Within both parts of the analysis, our proof relies on comparison arguments with other auxiliary chains.

2.1. *Reducing to two decks.* We use 2 auxiliary chains: a weighted and unweighted *deck-averaging* process. The use of a deck-averaging process is similar to the proof approach of Cancrini, Caputo and Martinelli [1] for the analysis of the so-called $L$-reversal chain.

The deck-averaging process is a continuous time Markov chain where the state space is again all the possible partitions of the $n$ cards into $k$ decks. For each pair of decks, there is an independent Poisson clock. In the unweighted process, all clocks have the same rate, whereas for the weighted process, the clock for decks $i$ and $j$ has rate $n_i + n_j$, where $n_\ell$ denotes the number of cards in deck $\ell$. The proof approach of [1] immediately yields that the unweighted deck-averaging process has a spectral gap of 1.

Our goal is to express the Dirichlet form and variance of the top-swap chain in terms of the sum of the 2-deck projections of the top-swap chain. It is straightforward to bound the Dirichlet form in terms of a sum over 2-deck Dirichlet forms, that is, where the configuration is fixed outside the two decks. Recall that the variance is independent of the chain. Thus, we can bound the variance by looking at the spectral gap of the unweighted deck-averaging process. This leaves us with the Dirichlet form for the unweighted deck-averaging process. We then use a nontrivial comparison argument to obtain the Dirichlet form for the weighted deck-averaging process. Finally, it is straightforward to bound the Dirichlet form of the weighted deck-averaging process as a sum over 2-deck variances. The result is a bound on the spectral gap of the ($k$-deck) top-swap chain in terms of the spectral gap of the 2-deck top-swap chain.

2.2. *Analysis of the* 2-*deck top-swap.* The basic idea is to compare the 2-deck top-swap with random transpositions. However, transpositions of two cards within the same deck are a problem. Roughly speaking, we can not efficiently "simulate" these transitions by top-swap transitions (see Remark 3.1 below for more details). Hence, we consider a transposition chain which only allows transpositions of cards in different decks. Moreover, if one of the cards is at the top of either deck, then instead of a transposition, the chain does the corresponding top-swap transition. It is straightforward to compare this modified transposition chain with the top-swap chain. To analyze this modified transposition chain, we compare to a final chain, which we call here the balanced-swap chain. The balanced-swap chain is also defined on two decks and has two types of transitions: swapping and rearranging. The swapping transition changes the size of the decks. In particular, we choose a random position in a deck, and move the cards above it, to the top of the other deck. The rearranging transition randomizes the card order while maintaining the current deck size. However, we only perform the second transition if the decks are balanced, that is, close to the same size. (This is made precise with the notion of $\delta$-balanced, for appropriately chosen $0 <$ $\delta < 1$.) It turns out that this later restriction is crucial for the comparison with the modified transposition chain.

Finally, we analyze the gap of the balanced-swap chain by adapting an analytical argument given by Cancrini, Martinelli, Roberto and Toninelli in their work on kinetically constrained models [2].

In the following we first describe the solution to the 2-deck problem in Section 3, then proceed with the deck-averaging process in Section 4, and finally in Section 5 we return to the analysis of the $k$-deck top-swap process. The matching upper bounds on the spectral gap are briefly discussed at the end of Section 5.

**3. Analysis of the 2-deck problem.** We consider two decks $I_1, I_2$ and we let $n_1, n_2$ denote the number of cards in each deck. Also, let $n = n_1 + n_2$ denote the total number of cards and let $\Omega_n$ be the set of all $(n+1)!$ arrangements of $n$ labeled cards in the two decks. We denote by $\eta$ a generic element of $\Omega_n$ (a configuration) and call $\mu$ the uniform probability measure on $\Omega_n$.

The 2-deck top-swap Markov chain is described as follows. At each step we choose two positions $r, s$, with each position drawn at random from the $n + 2$ available positions, $n_1 + 1$ from the first deck and $n_2 + 1$ from the second deck. As explained in the Introduction, if $r$ and $s$ belong to the same deck we do nothing. Otherwise, we swap the tops identified by positions $r, s$. The extra position added to each deck allows to swap an empty top. If the current configuration is $\eta$, we denote by $\mathcal{T}_{r,s}\eta$ the updated configuration. The Dirichlet form of the Markov chain can then be written as

$$(3.1) \qquad \mathcal{E}_2(f) = \frac{1}{2(n+2)^2} \sum_{r,s=1}^{n+2} \mu[(\mathcal{T}_{r,s}f - f)^2],$$

where $\mathcal{T}_{r,s}f(\eta) := f(\mathcal{T}_{r,s}\eta)$ for arbitrary functions $f : \Omega_n \to \mathbb{R}$, and $\mu[f]$ stands for expectation of a function $f$ w.r.t. the uniform probability $\mu$. The main result in this section is the following $O(n)$ estimate on the relaxation time $\tau(n,2)$ of the Markov chain described above.

THEOREM 3.1. *There exists $C > 0$ such that, for every $n \in \mathbb{N}$,*

$$(3.2) \qquad \tau(n,2) \leq Cn.$$

Recall that the estimate (3.2) is equivalent to showing that, for every function $f$,

$$(3.3) \qquad \text{Var}_\mu(f) \leq Cn\mathcal{E}_2(f).$$

We start by introducing a convenient notation. We add a card $*$ to mark the separation between the two decks and we call $x_* = x_*(\eta)$ the position of the $*$ in $\eta$; see Figure 1.

With this representation the top-swap transformations $\mathcal{T}_{r,s}$ can be rewritten by means of the transformations $\eta \to T_{i,j}\eta$, $i,j = 1, \ldots, n+1$, described below.

3.1. *Top-swap operators $T_{i,j}$.* If $i < x_* < j$, we call $T_{i,j}\eta$ the new configuration obtained by the ordinary top-swap. Namely, we collect the cards on top of $j$ from the second deck (including $j$) in a deck $D_j$ and the cards on top of $i$ from the first deck (including $i$) in a deck $D_i$. We then swap their positions, namely, $D_j$ goes above the position $(i-1)$ in the first deck and $D_i$ goes above the position $(j-1)$ in the second deck. Consider, for instance, the configuration $\eta = (3, 4, *, 5, 2, 1, 6)$ given in Figure 1 with $x_* = 3$. If we pick $i = 1$ and $j = 5$, we obtain $T_{1,5}\eta = (2, 1, 6, *, 5, 3, 4)$. If we choose $i = 2$ and $j = 7$, then $T_{2,7}\eta = (3, 6, *, 5, 2, 1, 4)$.

If $i = x_* < j$, then we move the deck $D_j$ above the position $(i-1)$ in the first deck. Therefore in our example, if $i = x_* = 3$ and, say, $j = 6$, then $T_{3,6}\eta = (3, 4, 1, 6, *, 5, 2)$.

Also, if $i < x_* = j$, we move the deck $D_i$ above the last card in the second deck. Thus, in our example if, say, $i = 1$ and $j = x_* = 3$, then $T_{1,3}\eta = (*, 5, 2, 1, 6, 3, 4)$.

Finally, if $i = j$ or $1 \leq i < j < x_*$ or $1 \leq x_* < i < j$, we do nothing, that is, we define $T_{i,j}\eta = \eta$ in these cases. To finish the definition of $T_{i,j}$ for all $i, j = 1, \ldots, n+1$, we set $T_{j,i} := T_{i,j}$ for every $i \leq j$.

Note that with these definitions we may rewrite the Dirichlet form (3.1) as

$$(3.4) \quad \mathcal{E}_2(f) = \frac{1}{2(n+2)^2} \sum_{i,j=1}^{n+1} \mu[(T_{i,j}f - f)^2],$$

where $T_{i,j}f(\eta) := f(T_{i,j}\eta)$. Here we use the fact that the transformations $T$ and $\mathcal{T}$ coincide when the chosen positions correspond to ordinary cards and are different if one of the tops to be swapped is empty. Note that the positions range from 1 to $n+2$ in (3.1) and from 1 to $n+1$ in (3.4). However, when one of the positions is the top of a deck in (3.1), then the second position must be in the other deck and there is no overcounting since in (3.4) we always perform the swap when one of the positions coincides with $x_*$.

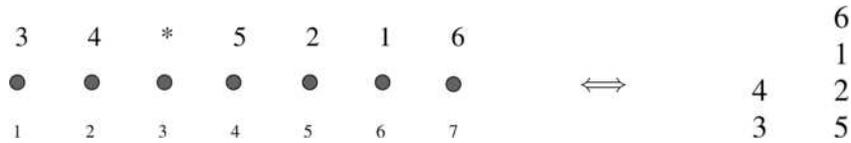

FIG. 1. *Two decks with $n = 6$. The configuration $\eta = (3, 4, *, 5, 2, 1, 6)$, with $n_1 = 2, n_2 = 4$ and $x_* = 3$.*

3.2. *Modified transpositions.* A first step in the proof of Theorem 3.1 is to compare the process described above to a random transposition-type Markov chain that is defined as follows. Let $E_{i,j}$ denote the ordinary transposition operator between positions $i$ and $j$. Namely, $(E_{i,j}\eta)_\ell = \eta_\ell$ for $i,j \neq \ell$ and $(E_{i,j}\eta)_i = \eta_j$, $(E_{i,j}\eta)_j = \eta_i$. Next, define the modified transposition $\widetilde{E}_{i,j}$, for $i < j$, by

$$(3.5) \qquad \widetilde{E}_{i,j} = \begin{cases} E_{i,j}, & \text{if } i < x_* < j, \\ T_{i,j}, & \text{if } i = x_* < j \text{ or } i < x_* = j, \\ 1, & \text{if } i < j < x_* \text{ or } x_* < i < j. \end{cases}$$

Also, set $\widetilde{E}_{i,i} = 1$ and $\widetilde{E}_{i,j} = \widetilde{E}_{j,i}$ if $i > j$. Here $T_{i,j}$ is the top-swap operation defined previously. In words, $\widetilde{E}_{i,j}$ is an ordinary transposition if the two chosen cards (in position $i$ and $j$ resp.) are in different decks. If they are in the same deck, nothing happens. Finally, if one of the cards is the $*$, then a top-swap transformation is performed. In analogy with (3.4), we define the Dirichlet form of the modified random transpositions by

$$(3.6) \qquad \mathcal{D}_2(f) = \frac{1}{2(n+2)^2} \sum_{i,j=1}^{n+1} \mu[(\widetilde{E}_{i,j}f - f)^2].$$

It turns out (see Lemma 3.4 below) that the comparison $\mathcal{D}_2 \leq C\mathcal{E}_2$ is easy. We then have to show that $\mathcal{D}_2$ has a spectral gap not smaller than constant times $n^{-1}$. This is accomplished by a further comparison with a new process with Dirichlet form $\mathcal{F}_\delta(f)$ defined below.

REMARK 3.1. We observe that comparison with ordinary random transpositions is not sufficient for our purpose. Let us call $\mathcal{E}_2^{\mathrm{RT}}(f)$ the Dirichlet form

$$\mathcal{E}_2^{\mathrm{RT}}(f) = \frac{1}{2(n+2)^2} \sum_{i,j=1}^{n+1} \mu[(E_{i,j}f - f)^2],$$

associated to pure transpositions $E_{i,j}$. Note that cards are not required to belong to different decks. The natural normalization here, as in (3.6), would be $2(n+1)^2$ instead of $2(n+2)^2$. However, the latter has been preferred to make the analogy with (3.4) more apparent. Then it is simple to check that taking, for example, $f(\eta) = 1$ if $\eta = (1, 2, \ldots, n, *)$ and $f(\eta) = 0$ otherwise, the ratio $\mathcal{E}_2^{\mathrm{RT}}(f)/\mathcal{E}_2(f)$ is of order $n$.

3.3. *The balanced-swap chain.* In words, this new process is described as follows. Let us fix $\delta \in (0, 1/2)$. We have two independent Poisson processes $\sigma$ and $\tau$ with mean 1. At the arrival times of $\sigma$ we choose uniformly at random a position $i$ and then update the current configuration $\eta$ by means of the transformation $\widetilde{E}_{i,x_*} = T_{i,x_*}$ [recall that $x_*(\eta)$ denotes the position of the $*$ in the configuration $\eta$]. At the arrival times of $\tau$ we look at the current value of $x_*(\eta)$ and, if

$$\delta \leq \frac{x_*(\eta)}{n+1} \leq (1-\delta), \tag{3.7}$$

then, conditioned on this value we put the system in equilibrium, that is, we choose the new configuration $\eta'$ uniformly at random among all the $n!$ configurations with the $*$ in the current position $x_*(\eta)$. If, instead, (3.7) does not hold, then we do nothing. We say that $\eta$ is $\delta$-balanced if (3.7) holds. We call $p_\delta$ the $\mu$-probability that this happens. Since $\mu$ is uniform, $p_\delta = 1 - 2\delta + O(1/n)$.

We write $\mathrm{Var}_\mu(f|x_*)$ for the variance of a function $f:\Omega_n \to \mathbb{R}$ w.r.t. the conditional probability $\mu(\cdot|x_*)$ obtained from $\mu$ by conditioning on the position of the $*$. The Dirichlet form $\mathcal{F}_\delta$ of the process described above is then given by

$$\mathcal{F}_\delta(f) = \mu[\chi \, \mathrm{Var}_\mu(f|x_*)] + \frac{1}{2(n+1)} \sum_{i=1}^{n+1} \mu[(T_{i,x_*}f - f)^2], \tag{3.8}$$

where

$$\chi(\eta) = \begin{cases} 1, & \text{if } \eta \text{ is } \delta\text{-balanced}, \\ 0, & \text{otherwise}. \end{cases} \tag{3.9}$$

We will show (see Lemma 3.3) that one has the comparison $\mathcal{F}_\delta \leq Cn\mathcal{D}_2$. The original problem (3.3) is then reduced to estimating from below the spectral gap of $\mathcal{F}_\delta$. To this end, one can use a coupling argument, but the following proof seems to give a better estimate.

LEMMA 3.2. *For all $n \in \mathbb{N}$, $f:\Omega_n \to \mathbb{R}$, we have*

$$\mathrm{Var}_\mu(f) \leq \gamma_\delta \mathcal{F}_\delta(f), \tag{3.10}$$

*where $\gamma_\delta := (1 - \sqrt{1-p_\delta})^{-1}$.*

PROOF. We have to prove that the spectral gap $\lambda_*$ of the Dirichlet form $\mathcal{F}_\delta(\cdot)$ satisfies

$$\lambda_* \geq 1 - \sqrt{1-p_\delta}. \tag{3.11}$$

The generator associated to $\mathcal{F}_\delta(\cdot)$ is the operator $\mathcal{L}$ acting on functions $f$ as

$$\mathcal{L}f = \chi(\mu[f|x_*] - f) + Af - f, \qquad Af := \frac{1}{n+1}\sum_{i=1}^{n+1} T_{i,x_*}f. \tag{3.12}$$

Note that $\mathcal{L}$ is self-adjoint and nonnegative in $L^2(\mu)$ and

$$-\mu[f\mathcal{L}f] = \mathcal{F}_\delta(f). \tag{3.13}$$

We consider the eigenvalues $\lambda$ and eigenfunctions $f_\lambda$ of $-\mathcal{L}$:

$$-\mathcal{L}f_\lambda = \lambda f_\lambda. \tag{3.14}$$

Here $\lambda = 0$ corresponds to $f_\lambda = $ constant and any $f_\lambda$ in (3.14) with $\lambda \neq 0$ satisfies $\mu[f_\lambda] = 0$. By definition, $\lambda_*$ is the smallest nonzero $\lambda$ such that (3.14) holds.

Observe that $\lambda = 1$ is an eigenvalue. Indeed, it suffices to consider $f(\eta) = \varphi(x_*(\eta))$ for an arbitrary $\varphi : \{1, \ldots, n+1\} \to \mathbb{R}$ such that $\frac{1}{n+1} \sum_{i=1}^{n+1} \varphi(i) = 0$. In this case, $f$ is a function of $x_*$ only and, therefore, $\mu[f|x_*] = f$ and the mean-zero property implies $Af = 0$, thus, (3.14) gives $\lambda = 1$.

We may assume without loss of generality that $\lambda_* < 1$, since otherwise there is nothing to prove in (3.11). We claim that $f_{\lambda_*}$ is such that

$$\mu[f_{\lambda_*}|x_*] = 0. \tag{3.15}$$

We first need the following property. For any function $g : \Omega_n \to \mathbb{R}$, we have

$$\mu[Ag|x_*] = \mu[g], \tag{3.16}$$

that is, $\mu[Ag|x_*]$ does not depend on the value of $x_*$. To prove (3.16), we write, for every $j$,

$$\mu[Ag|x_* = j] = \frac{1}{n+1} \sum_{i=1}^{n+1} \frac{\mu[(T_{i,j}g)1_{\{x_*=j\}}]}{\mathbb{P}_\mu[x_* = j]}.$$

Now we show that for any pair of positions $i, j$,

$$\mu[(T_{i,j}g)1_{\{x_*=j\}}] = \mu[g1_{\{x_*=i\}}]. \tag{3.17}$$

Indeed, if $\Omega_n^\ell := \{\eta \in \Omega_n : x_*(\eta) = \ell\}$, then it is not hard to see that $T_{i,j} : \Omega_n^j \to \Omega_n^i$ is a bijection: namely, if, for example, $i < j$, for every $\eta \in \Omega_n^j$ and for every $\xi \in \Omega_n^i$, we have $T_{i,n+i-j+2}\xi \in \Omega_n^j$, $T_{i,j}\eta \in \Omega_n^i$ and

$$\begin{aligned} T_{i,n+i-j+2}T_{i,j}\eta &= \eta, \\ T_{i,j}T_{i,n+i-j+2}\xi &= \xi. \end{aligned} \tag{3.18}$$

Therefore,

$$\mu[(T_{i,j}g)1_{\{x_*=j\}}] = \frac{1}{(n+1)!} \sum_{\eta \in \Omega_n^j} g(T_{i,j}\eta)$$

$$= \frac{1}{(n+1)!} \sum_{\xi \in \Omega_n^i} g(\xi) = \mu[g1_{\{x_*=i\}}].$$

This proves (3.17). Using $\mathbb{P}_\mu[x_* = j] = \frac{1}{n+1}$ and (3.17), we see that

$$\mu[Ag|x_* = j] = \sum_{i=1}^{n+1} \mu[g 1_{\{x_*=i\}}] = \mu[g],$$

which implies (3.16).

We turn to the proof of (3.15). Taking the $\mu[\cdot|x_*]$-expectation in (3.14), using (3.16) and $\mu[f_{\lambda_*}] = 0$, we see that

$$\lambda_* \mu[f_{\lambda_*}|x_*] = -\mu[\mathcal{L} f_{\lambda_*}|x_*]$$
$$= \mu[(f_{\lambda_*} - A f_{\lambda_*})|x_*] = \mu[f_{\lambda_*}|x_*].$$

Since $\lambda_* < 1$, this implies (3.15).

Using (3.15) in (3.14), we obtain

(3.19) $$A f_{\lambda_*} = (1 - \lambda_* + \chi) f_{\lambda_*} =: \psi f_{\lambda_*}.$$

Note that $\psi$ is a positive function of $x_*$ only. In particular, taking absolute values in (3.19), we have

(3.20) $$\mu[|A f_{\lambda_*}||x_*] = \psi(x_*) \mu[|f_{\lambda_*}||x_*].$$

Now, using $|A f_{\lambda_*}| \leq A|f_{\lambda_*}|$ and (3.16),

$$\mu[|A f_{\lambda_*}||x_*] \leq \mu[A|f_{\lambda_*}||x_*] = \mu[|f_{\lambda_*}|].$$

Therefore, dividing by $\psi$ and taking $\mu$-expectation in (3.20), we arrive at

(3.21) $$\mu[\psi^{-1}] \geq 1.$$

Inequality (3.21) can be written explicitly as

$$p_\delta \frac{1}{2 - \lambda_*} + (1 - p_\delta) \frac{1}{1 - \lambda_*} \geq 1.$$

It is easily seen that this is satisfied iff $\lambda_*^2 - 2\lambda_* + p_\delta \leq 0$, or $\lambda_* \geq 1 - \sqrt{1 - p_\delta}$ (since $\lambda_* < 1$ by assumption). This ends the proof of Lemma 3.2. □

REMARK 3.2. The result in Lemma 3.2 shows that if $p_\delta = 1$, that is, if there is no constraint and $\chi = 1$ (or $\delta = 0$), then the gap is equal to 1. It is, however, crucial for us that we are able to prove the spectral gap bound of Lemma 3.2 *with* the constraint $\chi < 1$, since we rely on this constraint to perform the comparison of Lemma 3.3 below.

3.4. *From balanced-swap moves to modified transpositions.* The next lemma shows how to "simulate" balanced-swap moves by means of modified transpositions.

LEMMA 3.3. *For each $\delta \in (0, 1/2)$, there exists $C_\delta > 0$ such that, for every $n$ and every function $f$ on $\Omega_n$,*

$$\mathcal{F}_\delta(f) \leq C_\delta n \mathcal{D}_2(f). \tag{3.22}$$

PROOF. Let us fix an arbitrary value of $x_*$. We want to estimate $\text{Var}(f|x_*)$ by using only terms of the form $\mu[(E_{i,\ell} - f)^2|x_*]$ with $i < x_* < \ell$.

Since $\mu(\cdot|x_*)$ is nothing but the uniform probability on $n!$ permutations (arrangements of $n$ cards in two decks of given size), recalling the spectral gap of the ordinary random transposition chain (see, e.g., [1]) we know that

$$\text{Var}_\mu(f|x_*) \leq \frac{1}{n} \sum_{i \neq x_*} \sum_{j \neq x_*} \mu[(E_{i,j}f - f)^2 | x_*]. \tag{3.23}$$

We now assume that $\delta n \leq x_* \leq (1-\delta)n$, for $\delta \in (0, 1/2)$ [see (3.7)]. Consider the case where $i, j$ are such that $i < x_*$ and $j < x_*$ (both cards are in the first deck). We observe that, for any $\ell > x_*$ we can write

$$E_{i,j}f - f = E_{i,\ell} E_{j,\ell} E_{i,\ell} f - f$$
$$= (E_{i,\ell}f_2 - f_2) + (E_{j,\ell}f_1 - f_1) + (E_{i,\ell}f - f),$$

where $f_1 := E_{i,\ell}f$ and $f_2 := E_{j,\ell}E_{i,\ell}f$. Therefore,

$$\mu[(E_{i,j}f - f)^2 | x_*]$$
$$\leq 3\mu[(E_{i,\ell}f_2 - f_2)^2 | x_*]$$
$$+ 3\mu[(E_{j,\ell}f_1 - f_1)^2 | x_*] + 3\mu[(E_{i,\ell}f - f)^2 | x_*].$$

From the invariance of $\mu(\cdot|x_*)$ under transpositions we see that

$$\mu[(E_{i,\ell}f_2 - f_2)^2 | x_*] = \mu[(E_{i,\ell}f - f)^2 | x_*],$$
$$\mu[(E_{j,\ell}f_1 - f_1)^2 | x_*] = \mu[(E_{j,\ell}f - f)^2 | x_*].$$

Therefore, averaging over all $n + 1 - x_*$ positions $\ell$ such that $\ell > x_*$, we obtain

$$\mu[(E_{i,j}f - f)^2 | x_*]$$
$$\leq \frac{1}{(n + 1 - x_*)} \sum_{\ell > x_*} \mu[6(E_{i,\ell}f - f)^2 + 3(E_{j,\ell}f - f)^2 | x_*].$$

Symmetrizing over $i, j$, the inequality above becomes

$$\mu[(E_{i,j}f - f)^2|x_*]$$
$$\leq \frac{9}{2(n+1-x_*)} \sum_{\ell > x_*} \mu[(E_{i,\ell}f - f)^2 + (E_{j,\ell}f - f)^2|x_*].$$

Summing over all pairs $i, j$ with $i < x_*$ and $j < x_*$, we obtain

$$\sum_{i<x_*} \sum_{j<x_*} \mu[(E_{i,j}f - f)^2|x_*]$$

(3.24)
$$\leq \frac{9(x_* - 1)}{(n+1-x_*)} \sum_{i<x_*} \sum_{\ell>x_*} \mu[(E_{i,\ell}f - f)^2|x_*],$$

since there are $x_* - 1$ positions $j$ such that $j < x_*$.

Repeating the reasoning leading to (3.24) for pairs $i, j$ with $i > x_*$ and $j > x_*$, we have

(3.25)
$$\sum_{i>x_*} \sum_{j>x_*} \mu[(E_{i,j}f - f)^2|x_*]$$
$$\leq \frac{9(n+1-x_*)}{(x_* - 1)} \sum_{i>x_*} \sum_{\ell<x_*} \mu[(E_{i,\ell}f - f)^2|x_*].$$

Inserting (3.24) and (3.25) in (3.23), we obtain

$$\operatorname{Var}(f|x_*) \leq \frac{2}{n} \sum_{i<x_*} \sum_{j>x_*} \mu[(E_{i,j}f - f)^2|x_*]$$
$$+ \frac{9}{n} \left\{ \frac{(n+1-x_*)}{(x_*-1)} + \frac{(x_*-1)}{(n+1-x_*)} \right\} \sum_{i<x_*} \sum_{j>x_*} \mu[(E_{i,j}f - f)^2|x_*]$$
$$\leq \frac{1}{n} \left\{ 2 + \frac{9n^2}{(n+1-x_*)(x_*-1)} \right\} \sum_{i<x_*} \sum_{j>x_*} \mu[(E_{i,j}f - f)^2|x_*],$$

where we used $(n-t)/t + t/(n-t) \leq n^2/t(n-t)$, for $t \in (0, n)$. Note that if $x_*$ satisfies (3.7), then $(n+1-x_*)(x_*-1) \geq \delta(1-\delta)(n+1)^2$. In conclusion, recalling that $\chi$ forces $x_*$ to satisfy (3.7), we see that, for any value of $x_*$,

$$\chi \operatorname{Var}(f|x_*) \leq \frac{1}{n} \left\{ 2 + \frac{9n^2}{\delta(1-\delta)(n+1)^2} \right\} \sum_{i<x_*} \sum_{j>x_*} \mu[(E_{i,j}f - f)^2|x_*].$$

Taking $\mu$-expectation and recalling the definition (3.8) of $\mathcal{F}_\delta(f)$, we see that

(3.26) $\quad \mathcal{F}_\delta(f) \leq \left\{ 2 + \frac{9n^2}{\delta(1-\delta)(n+1)^2} \right\} \frac{(n+2)^2}{n} \mathcal{D}_2(f) =: C_\delta n \mathcal{D}_2(f).$

Note that $C_\delta \leq 2 + 9/[\delta(1-\delta)] + O(1/n)$. □

3.5. *From modified transpositions to top-swaps.* Here we use the path-comparison technique of Diaconis and Saloff-Coste [4] to estimate the Dirichlet form $\mathcal{D}_2$ in terms of the Dirichlet form $\mathcal{E}_2$; see (3.4) and (3.6) for their definitions. The next lemma is based on the rather obvious fact that an ordinary transposition can be "simulated" by two top-swap moves.

LEMMA 3.4. *For any function $f$,*

$$\mathcal{D}_2(f) \leq 5\mathcal{E}_2(f). \tag{3.27}$$

PROOF. Let $\eta$ be a given configuration of cards and let $i, j$ be two given positions such that $i < j$. If $i < x_*(\eta) = j$, then $\widetilde{E}_{i,j} = T_{i,j}$, so we only need to take care of the case $i < x_*(\eta) < j$. In this case $\widetilde{E}_{i,j}$ is an ordinary transposition $E_{i,j}$ and this can be written by way of two top-swap operations. Namely, suppose that $\eta \in \Omega_n^\ell$, that is, $\eta \in \Omega_n$ is such that $x_*(\eta) = \ell$. Then it is not hard to see that

$$E_{i,j}\eta = T_{i+1,n+i-\ell+3}T_{i,j}\eta. \tag{3.28}$$

Note that $n + i - \ell + 3 \leq n + 2$ since we are assuming $i < x_*(\eta) = \ell$. We adopt the convention that if $n + i - \ell + 3 = n + 2$, then $T_{i+1,n+i-\ell+3}\xi$ changes the configuration $\xi$ by moving the top above (and including) $i + 1$ from deck 1 to the top of deck 2.

The identity (3.28) is easily understood with the help of a picture. For instance, consider the configuration $\eta = (3, 4, *, 5, 2, 1, 6)$ of Figure 1. $E_{1,6}\eta$ (i.e., the transposition of cards labeled 3 and 1) is obtained by first applying $T_{1,6}$ to obtain $T_{1,6}\eta = (1, 6, *, 5, 2, 3, 4)$ and then applying $T_{2,7}$. To see a case where $n + i - \ell + 3 = n + 2$, consider, for example, $E_{2,5}\eta$ (the transposition of cards labeled 4 and 2). This is obtained by doing first $T_{2,5}\eta = (3, 2, 1, 6, *, 5, 4)$ and then applying $T_{3,8}$ to this last configuration, since, by our convention, $T_{3,8}$ is the same as $T_{3,5}$.

Using (3.28), we have, for all $i < \ell < j$,

$$(E_{i,j}f(\eta) - f(\eta))^2 1_{\{x_*=\ell\}} \leq 2(T_{i+1,n+i-\ell+3}f_1 - f_1)^2 1_{\{x_*=\ell\}} + 2(T_{i,j}f - f)^2 1_{\{x_*=\ell\}}, \tag{3.29}$$

where $f_1(\eta) := f(T_{i,j}\eta)$ for every $\eta \in \Omega_n$. An argument similar to that used in the proof of (3.17) gives

$$\mu[(T_{i+1,n+i-\ell+3}f_1 - f_1)^2 1_{\{x_*=\ell\}}] = \mu[(T_{i+1,n+i-\ell+3}f - f)^2 1_{\{x_*=n+i-j+2\}}]. \tag{3.30}$$

From (3.30) it follows that

$$\sum_{i,j:i<j} \sum_{\ell=i+1}^{j-1} \mu[(T_{i+1,n+i-\ell+3}f_1 - f_1)^2 1_{\{x_*=\ell\}}]$$

(3.31)
$$= \sum_{i,j:i<j} \sum_{\ell=i+1}^{j-1} \mu[(T_{i+1,n+i-\ell+3}f - f)^2 1_{\{x_*=n+i-j+2\}}]$$

$$\leq \sum_{i,j:i<j} \mu[(T_{i,j}f - f)^2].$$

Inserting (3.31) in (3.29), we have shown that

(3.32) $$\sum_{i,j=1}^{n+1} \mu[(\widetilde{E}_{i,j}f - f)^2 1_{\{i<x_*<j\}}] \leq 4 \sum_{i,j:i<j} \mu[(T_{i,j}f - f)^2].$$

The identity

$$\mathcal{D}_2(f) = \frac{1}{(n+2)^2} \sum_{i,j=1}^{n+1} \mu[(\widetilde{E}_{i,j}f - f)^2 1_{\{i<x_*<j\}}]$$

$$+ \frac{1}{(n+2)^2} \sum_{i=1}^{n+1} \mu[(T_{i,x_*}f - f)^2]$$

and (3.32) then imply

$$\mathcal{D}_2(f) \leq \frac{2}{(n+2)^2} \sum_{i,j=1}^{n+1} \mu[(T_{i,j}f - f)^2]$$

$$+ \frac{1}{(n+2)^2} \sum_{i}^{n+1} \mu[(T_{i,x_*}f - f)^2] \leq 5\mathcal{E}_2(f). \quad \square$$

3.6. *Proof of Theorem* 3.1 *completed.* As a consequence of Lemma 3.2, Lemma 3.3 and Lemma 3.4, we know that, for all $n$ and all $f:\Omega_n \to \mathbb{R}$,

(3.33) $$\mathrm{Var}_\mu(f) \leq 5C_\delta \gamma_\delta n \mathcal{E}_2(f).$$

We can still choose the value of $\delta \in (0, 1/2)$ ($C_\delta$ is minimized at $\delta = 1/2$, while $\gamma_\delta$ is minimized at $\delta = 0$).

Note that for $\delta = 0.25$ one has $p_\delta = 0.5 + O(1/n)$, $\gamma_\delta = \frac{\sqrt{2}}{\sqrt{2}-1} + O(1/n)$ and $C_\delta \leq 50 + O(1/n)$. Therefore, (3.33) gives $\mathrm{Var}_\mu(f) \leq Cn\mathcal{E}_2(f)$ with $C \leq 875 + O(1/n)$. This ends the proof of Theorem 3.1.

3.7. *Random transpositions with constraint.* The analysis above can be used to obtain a spectral gap estimate for a random transposition model with the constraint that two cards need to be in two different decks to be transposed. This is the process with Dirichlet form $\mathcal{D}_2^T(f)$ as defined in (3.6), with the only difference that here the transformations $\widetilde{E}_{i,j}$ are given by

$$(3.34) \qquad \widetilde{E}_{i,j}^T = \begin{cases} E_{i,j}, & \text{if } i \leq x_* \leq j, \\ 1, & \text{if } i < j < x_* \text{ or } x_* < i < j, \end{cases}$$

with $\widetilde{E}_{i,i}^T = 1$ and $\widetilde{E}_{i,j}^T = \widetilde{E}_{j,i}^T$ if $i > j$. Namely, when $i = x_*$, we simply transpose rather than doing the top-swap as in (3.5).

THEOREM 3.5. *There exists $C < \infty$ such that, for every $n \in \mathbb{N}$, for every function $f: \Omega_n \to \mathbb{R}$, we have*

$$(3.35) \qquad \mathrm{Var}_\mu(f) \leq Cn\mathcal{D}_2^T(f).$$

PROOF. We repeat the same arguments used in the proof of Theorem 3.1. We only need to modify the definition of the Dirichlet form $\mathcal{F}_\delta(f)$ which we replace here with

$$(3.36) \qquad \mathcal{F}_\delta^T(f) = \mu[\chi \mathrm{Var}_\mu(f|x_*)] + \frac{1}{2(n+1)} \sum_{i=1}^{n+1} \mu[(E_{i,x_*}f - f)^2].$$

It is not difficult to check that all our arguments in Lemma 3.2 apply with no modifications to this case. Moreover, the same applies to Lemma 3.3. Note that here we do not need the extra comparison of Lemma 3.4. In particular, one has that the constant $C$ in (3.35) satisfies $C \leq 175 + O(1/n)$. □

**4. Deck-averages.** We consider the following setting. There are $k$ decks $I_1, \ldots, I_k$ and a total of $n$ cards. The state space $\Omega$ is the set of all

$$\frac{(n+k-1)!}{(k-1)!}$$

arrangements of $n$ cards in $k$ decks. $\eta$ will denote the general random element of $\Omega$ and $\nu$ will be the uniform probability over $\Omega$. Also, we use $\pi_i \eta$ to denote the projection of the configuration $\eta$ onto the configuration of cards in deck $I_i$.

4.1. *Deck-averaging process.* The deck-averaging process is the continuous time Markov chain obtained as follows. At each deck independently there is a rate-1 Poisson clock. When deck $I_i$ rings, we choose uniformly another deck $I_j$ and the cards in the two decks $I_i \cup I_j$ are rearranged uniformly at random (among the two decks). The Dirichlet form of this process is

$$(4.1) \qquad \bar{\mathcal{D}}(f) = \frac{1}{k} \sum_{i=1}^{k} \sum_{j=1}^{k} \nu[(A_{i,j}f)^2],$$

where the averaging gradient is defined by

$$
\begin{aligned}
A_{i,j}f(\sigma) &= \sum_{\xi \in \Omega}(f(\xi) - f(\sigma))\nu[\eta = \xi | \pi_\ell \eta = \pi_\ell \sigma, \ell \neq i, j] \\
&= \nu[f | \pi_\ell \eta, \ell \neq i, j](\sigma) - f(\sigma) \\
&=: \nu_{i,j}[f](\sigma) - f(\sigma),
\end{aligned}
\tag{4.2}
$$

for $i \neq j$ and $\sigma \in \Omega$. We set $A_{i,i}f = 0$. Here $\nu_{i,j} = \nu[\cdot | \pi_\ell \eta, \ell \neq i, j]$ is the conditional probability obtained by freezing the configuration in all decks $I_\ell$, $\ell \neq i, j$. Also, note that

$$
\nu[(A_{i,j}f)^2] = \nu[\mathrm{Var}_{i,j}(f)], \tag{4.3}
$$

where $\mathrm{Var}_{i,j}$ stands for the variance with respect to the probability $\nu_{i,j}$, with the convention that $\mathrm{Var}_{i,i}(f) = 0$.

The result of Proposition 2.3 in [1] is that the spectral gap of the deck-average Markov chain equals 1. In particular, for every $f : \Omega \to \mathbb{R}$, we have

$$
\mathrm{Var}_\nu(f) \leq \bar{\mathcal{D}}(f), \tag{4.4}
$$

with $\mathrm{Var}_\nu(f)$ denoting the variance of $f$ w.r.t. $\nu$. Strictly speaking, the proof of Proposition 2.3 in [1] does not apply here because the sizes of the decks are not fixed in our model. However, it is not hard to adapt that proof to show that the result (4.4) holds. This only requires small modifications, the point being that the spectrum of the operator $\mathcal{K}$ used there satisfies the right bounds in our setting; see Lemma 4.13 below and the remark following it for more details.

4.2. *Weighted deck-averaging process.* The weighted version of the deck-averaging chain is defined by the Dirichlet form

$$
\mathcal{D}(f) = \frac{1}{k} \sum_{i=1}^{k} \sum_{j=1}^{k} \nu[(n_i + n_j)\mathrm{Var}_{i,j}(f)], \tag{4.5}
$$

where $n_i = n_i(\eta)$ stands for the number of cards in deck $i$. A possible interpretation for (4.5) is that each *card* is equipped with a rate-1 Poisson clock; when a card $c$ rings a deck, $I_j$ is chosen uniformly at random; if $c \in I_j$, then nothing happens, while if $c \notin I_j$, that is, $c \in I_i$ for some $i \neq j$, then all the cards in $I_i \cup I_j$ are rearranged uniformly at random. Apart from a factor 2, the Dirichlet form of this process is given by (4.5). The main result here is the following spectral gap estimate.

PROPOSITION 4.1. *For every* $f : \Omega \to \mathbb{R}$,

$$
\mathrm{Var}_\nu(f) \leq \frac{6k}{n}\mathcal{D}(f). \tag{4.6}
$$

PROOF. Since $\sum_{i=1}^{k} n_i = n$, (4.4) and (4.3) imply

$$\text{Var}_\nu(f) \leq \frac{1}{kn} \sum_{i,j=1}^{k} \sum_{\ell=1}^{k} \nu[n_\ell \text{Var}_{i,j}(f)]. \tag{4.7}$$

We claim that, for every fixed triple $i, j, \ell$,

$$\begin{aligned}
\tfrac{1}{2}\nu[n_\ell \text{Var}_{i,j}(f)] &\leq \nu[(n_i + n_j)\text{Var}_{i,j}(f)] \\
&\quad + \nu[(n_i + n_\ell)\text{Var}_{i,\ell}(f)] \\
&\quad + \nu[(n_j + n_\ell)\text{Var}_{j,\ell}(f)].
\end{aligned} \tag{4.8}$$

Once we have (4.8), from the estimate in (4.7), we obtain

$$\text{Var}_\nu(f) \leq \frac{6}{n} \sum_{i,j=1}^{k} \nu[(n_i + n_j)\text{Var}_{i,j}(f)],$$

which is the same as (4.6).

We turn to the proof of the claim (4.8). We can assume that $i, j, \ell$ are three distinct labels (i.e., $\ell \neq i$, $\ell \neq j$, $i \neq j$) since the statement is obviously true otherwise. Let $m = n_i + n_j + n_\ell$ and let $\mu_\Delta$ denote the conditional probability $\nu[\cdot | \pi_v \eta, v \neq i, j, \ell]$ obtained by freezing all the decks $I_v$, $v \neq i, j, \ell$. We are going to prove that

$$\begin{aligned}
\frac{m}{2} \text{Var}_{\mu_\Delta}(f) &\leq \mu_\Delta[(n_i + n_j)\text{Var}_{i,j}(f)] \\
&\quad + \mu_\Delta[(n_i + n_\ell)\text{Var}_{i,\ell}(f)] \\
&\quad + \mu_\Delta[(n_j + n_\ell)\text{Var}_{j,\ell}(f)].
\end{aligned} \tag{4.9}$$

Here $\text{Var}_{\mu_\Delta}(f)$ is the variance of $f$ w.r.t. $\mu_\Delta$. This depends on the given configuration $\eta$ through the projections $\pi_v \eta$ for $v \neq i, j, \ell$. Note that once the variables $\pi_v \eta$, $v \neq i, j, \ell$ are given, $m$ becomes a constant, that is, $m$ is $\mu_\Delta$-a.s. constant.

Let us first show that (4.9) implies the claim (4.8). Take expectation w.r.t. $\nu$ in (4.9) and use the fact that $\nu[h] = \nu[\mu_\Delta[h]]$ for any function $h : \Omega \to \mathbb{R}$:

$$\begin{aligned}
\tfrac{1}{2}\nu[m \text{Var}_{\mu_\Delta}(f)] &\leq \nu[(n_i + n_j)\text{Var}_{i,j}(f)] \\
&\quad + \nu[(n_i + n_\ell)\text{Var}_{i,\ell}(f)] \\
&\quad + \nu[(n_j + n_\ell)\text{Var}_{j,\ell}(f)].
\end{aligned} \tag{4.10}$$

Observe that

$$\begin{aligned}
\text{Var}_{\mu_\Delta}(f) &= \mu_\Delta[f(f - \mu_\Delta[f])] \\
&= \mu_\Delta[f(f - \mu_\Delta[f|\pi_\ell \eta])] + \mu_\Delta[f(\mu_\Delta[f|\pi_\ell \eta] - \mu_\Delta[f])] \\
&= \mu_\Delta[\text{Var}_{i,j}(f)] + \text{Var}_{\mu_\Delta}[\mu_\Delta[f|\pi_\ell \eta]] \\
&\geq \mu_\Delta[\text{Var}_{i,j}(f)],
\end{aligned}$$

where we use the fact that $\mu_\Delta[f|\pi_\ell \eta]$ coincides with $\nu_{i,j}$.

We can now use the obvious bound $m \geq n_\ell$ to obtain

$$\nu[m \operatorname{Var}_{\mu_\Delta}(f)] \geq \nu[m \mu_\Delta[\operatorname{Var}_{i,j}(f)]]$$
$$= \nu[m \operatorname{Var}_{i,j}(f)] \geq \nu[n_\ell \operatorname{Var}_{i,j}(f)].$$

With (4.10), this implies (4.8).

It remains to prove (4.9). We introduce the projection $P_i$ defined by

$$P_i f = \mu_\Delta[f|\pi_i \eta].$$

Projectors $P_j, P_\ell$ are defined in a similar way. Observe that with these definitions, since $\mu_\Delta[\cdot|\pi_\ell \eta] = \nu_{i,j}$, we have

(4.11) $$\mu_\Delta[\operatorname{Var}_{i,j}(f)] = \mu_\Delta[f(1 - P_\ell)f].$$

Moreover, setting

$$\Gamma := (m - n_i)(1 - P_i) + (m - n_j)(1 - P_j) + (m - n_\ell)(1 - P_\ell),$$

we see that the r.h.s. of display (4.9) coincides with the quadratic form $\mu_\Delta[f\Gamma f]$.

We call $P$ the average:

(4.12) $$P = \tfrac{1}{3}(P_i + P_j + P_\ell).$$

Lemma 4.2 below shows that

(4.13) $$\mu_\Delta[fPf] \leq \tfrac{1}{2}\mu_\Delta[f^2],$$

for any function $f$ such that $\mu_\Delta[f] = 0$. For such $f$, the result (4.13) implies

$$\mu_\Delta[f\Gamma f] = 2m\mu_\Delta[f^2] - 3m\mu_\Delta[fPf] + \mu_\Delta[n_i f P_i f]$$
$$+ \mu_\Delta[n_j f P_j f] + \mu_\Delta[n_\ell f P_\ell f]$$
$$\geq \frac{m}{2}\mu_\Delta[f^2],$$

where we use

$$\mu_\Delta[n_i f P_i f] = \mu_\Delta[n_i (P_i f)^2] \geq 0.$$

The claim (4.9) follows since, by subtracting the mean $\mu_\Delta[f]$, we may restrict to mean zero functions to obtain

$$\operatorname{Var}_{\mu_\Delta}(f) \leq \frac{2}{m}\mu_\Delta[f\Gamma f]. \qquad \square$$

LEMMA 4.2. *Let $P$ denote the operator defined in (4.12). For every function $f$ such that $\mu_\Delta[f] = 0$, we have $\mu_\Delta[fPf] \leq \tfrac{1}{2}\mu_\Delta[f^2]$.*

PROOF. In the case of three decks, the proof of Proposition 2.2 of [1] [see display (27) and display (32) there] shows that

(4.14) $$\mu_\Delta[f(1-P)f] \geq \min\{\tfrac{1}{2}, \tfrac{2}{3}(1-\lambda_2)\}\mu_\Delta[f^2],$$

where $\lambda_2$ is the largest positive eigenvalue (other than $\lambda_0 = 1$; see below for more details) of the stochastic matrix $\mathcal{K}(\alpha,\beta)$ defined by

$$\mathcal{K}(\alpha,\beta) = \mu_\Delta[\pi_i\eta = \beta | \pi_j\eta = \alpha].$$

That is, $\mathcal{K}(\alpha,\beta)$ is the $\mu_\Delta$-conditional probability of having the configuration $\beta$ in $I_i$, given the configuration $\alpha$ in $I_j$.

The spectrum of $\mathcal{K}$ can be computed exactly as in Lemma 2.1 of [1]. The only difference here is that the size of the decks is not fixed. This does not pose any difficulty and we can proceed in the very same way. We have to compute the action of $\mathcal{K}$ on functions of the kind $\chi_{r_1}\cdots\chi_{r_n}$, where $\chi_r$ stands for the indicator function of the event {the card labeled $r$ belongs to deck $i$}. For instance, a simple computation shows that for one card $r$ one has

$$\mathcal{K}\chi_r = \tfrac{1}{2}(1-\chi_r).$$

Therefore, $\chi_r - \tfrac{1}{3}$ is an eigenfunction with eigenvalue $\lambda_1 = -\tfrac{1}{2}$. Similarly, for two distinct cards $r_1, r_2$, one can compute

$$\mathcal{K}\chi_{r_1}\chi_{r_2} = \tfrac{1}{4}(1-\chi_{r_1})(1-\chi_{r_2}).$$

Following the proof of Lemma 2.1 in [1], we then obtain that the spectrum of $\mathcal{K}$ consists of the eigenvalues $\lambda_\ell = (-\tfrac{1}{2})^\ell$, for $\ell = 0, 1, \ldots, m$, where $m$ is the total number of cards in the three decks. In particular, $\lambda_2 = \tfrac{1}{4}$ so that (4.14) implies $\mu_\Delta[f(1-P)f] \geq \tfrac{1}{2}\mu_\Delta[f^2]$ or, equivalently,

$$\mu_\Delta[fPf] \leq \tfrac{1}{2}\mu_\Delta[f^2]. \qquad \square$$

REMARK 4.1. More generally, for $k$ decks and $n$ cards, one obtains that the spectrum of $\mathcal{K}$ is given by

$$\lambda_\ell = \left(-\frac{1}{k-1}\right)^\ell,$$

$\ell = 0, 1, \ldots, n$.

**5. Top-swap process with $k$ decks.** The setting is as in the previous section. In analogy with the 2-deck process defined in (3.1), we consider here the Dirichlet form

(5.1) $$\mathcal{E}_k(f) = \frac{1}{2(n+k)^2} \sum_{r,s=1}^{n+k} \nu[(\mathcal{T}_{r,s}f - f)^2],$$

where we sum over all positions $r, s$ between 1 and $n + k$, since we added to each deck $I_j$ an extra position to allow the swap of an empty top. We shall call $\bar{I}_j$ the set of positions defined by the deck $I_j$ and its extra position.

The main result is that the relaxation time of this chain is $O(n+k)$:

THEOREM 5.1. *There exists a constant $C > 0$ such that, for every $n, k$, every $f$*

$$\operatorname{Var}_\nu(f) \leq C(n+k)\mathcal{E}_k(f). \tag{5.2}$$

PROOF. We start by decomposing $\mathcal{E}_k$ into 2-deck terms:

$$\mathcal{E}_k(f) = \frac{1}{4(n+k)^2} \sum_{i,j=1}^{k} \nu \left[ \sum_{r,s \in \bar{I}_i \cup \bar{I}_j} (\mathcal{T}_{r,s}f - f)^2 \right]. \tag{5.3}$$

If $\nu_{i,j}$ denotes the conditional probability in (4.2) (obtained by freezing all the decks $I_\ell$, $\ell \neq i, j$), then

$$\begin{aligned}
\mathcal{E}_k(f) &= \frac{1}{4(n+k)^2} \sum_{i,j=1}^{k} \nu \left[ \sum_{r,s \in \bar{I}_i \cup \bar{I}_j} \nu_{i,j}[(\mathcal{T}_{r,s}f - f)^2] \right] \\
&= \frac{1}{2(n+k)^2} \sum_{i,j=1}^{k} \nu[(n_i + n_j + 2)^2 \mathcal{E}_{i,j}(f)],
\end{aligned} \tag{5.4}$$

where $\mathcal{E}_{i,j}(f)$ is the Dirichlet form of the top-swap process on the two decks $i, j$ only:

$$\mathcal{E}_{i,j}(f) = \frac{1}{2(n_i + n_j + 2)^2} \nu \left[ \sum_{r,s \in \bar{I}_i \cup \bar{I}_j} \nu_{i,j}[(\mathcal{T}_{r,s}f - f)^2] \right].$$

From Theorem 3.1, we know that for all $f : \Omega \to \mathbb{R}$, for any pair of distinct decks $I_i, I_j$ with fixed total number of cards $n_i + n_j$,

$$\mathcal{E}_{i,j}(f) \geq \frac{1}{C(n_i + n_j)} \operatorname{Var}_{i,j}(f). \tag{5.5}$$

We can then use (5.4) to obtain

$$\mathcal{E}_k(f) \geq \frac{1}{2C(n+k)^2} \sum_{i,j=1}^{k} \nu[(n_i + n_j + 2) \operatorname{Var}_{i,j}(f)]. \tag{5.6}$$

Recalling the definition (4.5) of the weighted deck-averaging Dirichlet form $\mathcal{D}(f)$ we see that

$$\mathcal{E}_k(f) \geq \frac{1}{2C(n+k)^2} \left[ k\mathcal{D}(f) + 2 \sum_{i,j=1}^{k} \nu(\operatorname{Var}_{i,j}(f)) \right]. \tag{5.7}$$

Now using (4.3), (4.1) and (4.4), we get that

$$(5.8) \qquad \mathcal{E}_k(f) \geq \frac{k}{2C(n+k)^2}[\mathcal{D}(f) + 2\operatorname{Var}_\nu(f)].$$

Proposition 4.1 then implies

$$\mathcal{E}_k(f) \geq \frac{k}{2C(n+k)^2}\left[\frac{n}{6k}\operatorname{Var}_\nu(f) + 2\operatorname{Var}_\nu(f)\right]$$

$$= \frac{n+12k}{12C(n+k)^2}\operatorname{Var}_\nu(f)$$

$$\geq \frac{1}{12C(n+k)}\operatorname{Var}_\nu(f).$$

This concludes the proof. □

5.1. *k-deck random transpositions with constraints.* Here we consider the extension to $k$ decks of the problem discussed in Theorem 3.5. Namely, we define

$$(5.9) \qquad \mathcal{E}_k^T(f) = \frac{1}{2(n+k)^2}\sum_{r,s=1}^{n+k}\nu[(E_{r,s}f - f)^2],$$

where the transformations $E_{r,s}$ are interpreted as ordinary transpositions provided $r$ and $s$ belong to different extended decks $\bar{I}_i, \bar{I}_j$. If $r, s \in \bar{I}_i$ for some $i$, then $E_{r,s} = 1$.

The computations in the proof given above can be repeated step by step. Using Theorem 3.5, we can then obtain the following analogue of Theorem 5.1.

THEOREM 5.2. *There exists $C > 0$ such that for every $n, k$ and every $f$*

$$(5.10) \qquad \operatorname{Var}_\nu(f) \leq C(n+k)\mathcal{E}_k^T(f).$$

5.2. *Matching upper bounds on spectral gap.* Consider the $k$-deck top-swap Dirichlet form given by (5.1). We now show that, up to a constant factor, Theorem 5.1 is tight.

For $\eta \in \Omega$, let $n_1(\eta)$ denote as usual the number of cards in the first deck. Consider the indicator (test) function

$$(5.11) \qquad f(\eta) = \begin{cases} 1, & \text{if } n_1(\eta) = 0, \\ 0, & \text{if } n_1(\eta) > 0. \end{cases}$$

Then it is easy to check that $\nu[f] = (k-1)/(n+k-1)$, and that

$$\operatorname{Var}_\nu(f) = n(k-1)/(n+k-1)^2.$$

On the other hand,

$$\mathcal{E}_k(f) = \frac{1}{2(n+k)^2}\left[\frac{2(k-1)n}{(n+k-1)}\right]$$
$$= \frac{n(k-1)}{(n+k-1)(n+k)^2},$$
(5.12)

thus, showing that the spectral gap is at most

$$(n+k-1)/(n+k)^2 \leq 1/(n+k).$$

Note that in the first equality above we used the fact that, for $\eta$ with $n_1(\eta) = 0$, the contribution to the Dirichlet form comes from choosing $r = 1$ (viz., the position of the first marker) and $s$ to be any of the $n$ card positions (viz., a nonmarker position).

Finally, note that the same test function (with an identical computation) shows that the spectral gap of the $k$-deck random transpositions with constraints, defined using the Dirichlet form (5.9), is also at most $1/(n+k)$.

**Acknowledgment.** We thank Fabio Martinelli for pointing out the connection with [2] and for suggesting the argument of Lemma 3.2.

N. Bhatnagar
E. Vigoda
College of Computing
Georgia Institute of Technology
Atlanta, Georgia 30332
USA
E-mail: nand@cc.gatech.edu
vigoda@cc.gatech.edu

P. Caputo
Dipartimento di Matematica
Università di Roma Tre
Largo S. Murialdo 1
00146 Rome
Italy
E-mail: caputo@mat.uniroma3.it

P. Tetali
School of Mathematics
  and College of Computing
Georgia Institute of Technology
Atlanta, Georgia 30332
USA
E-mail: tetali@math.gatech.edu